\documentclass{easychair}

\usepackage{lineno}
\usepackage{doc}
\usepackage{booktabs}
\usepackage{amsmath,amsfonts}
\usepackage{algorithm}
\usepackage{algpseudocode}
\usepackage{graphicx}
\usepackage{multirow}
\usepackage{multicol}
\usepackage{amssymb}
\usepackage{subcaption}
\usepackage{xcolor}
\usepackage{diagbox}
\usepackage{enumitem}
\setlist[itemize]{leftmargin=*}
\usepackage{url}
\usepackage{threeparttable}
\usepackage{pifont}
\usepackage{nopageno}
\usepackage{bm}
\usepackage{scalerel}
\usepackage{makecell}
\usepackage[colorinlistoftodos]{todonotes}

\newcommand{\R}{\mathbb{R}}

\title{
Reduction of large-scale RLCk models \\ via low-rank balanced truncation
}

\author{
Christos Giamouzis,
Dimitrios Garyfallou,
Anastasis Vagenas,\\
and Nestor Evmorfopoulos
}

\institute{
 Dept. of Electrical and Computer Engineering, 
University of Thessaly, Volos, Greece \\
  \email{\{cgiamouzis, digaryfa, avagenas, nestevmo\}@e-ce.uth.gr}
}

\begin{document}
\maketitle

\begin{abstract}
Model order reduction (MOR) is an important step in the design process of integrated circuits. Specifically, the electromagnetic models extracted from modern complex designs result in a large number of passive elements that introduce limitations in the simulation process. MOR techniques based on balanced truncation (BT) can overcome these limitations by producing compact reduced-order models (ROMs) that approximate the behavior of the original models at the input/output ports. In this paper, we present a low-rank BT method that exploits the extended Krylov subspace and efficient implementation techniques for the reduction of large-scale models. Experimental evaluation on a diverse set of analog and mixed-signal circuits with millions of elements indicates that up to $\times$5.5 smaller ROMs can be produced with similar accuracy to ANSYS~RaptorX\texttrademark\ ROMs.
\end{abstract}
\vspace{-6.7pt}
\section{Introduction}
\label{sec:intro}

Electromagnetic model extraction plays a key role in the design and analysis of integrated circuits. The extracted models are simulated to accurately predict the behavior of the passive elements of the design. Model order reduction (MOR) can reduce the complexity of RLCk models with many elements ($>$1M) and ports ($>$10), while retaining an accurate approximation of the input and output behavior of the circuit~\cite{Odabasioglu1998, Antoulas2004}. Therefore, the simulation time of complex systems can be radically decreased by constructing reduced-order models (ROMs) of smaller dimensions that preserve the essential characteristics of the original models.

MOR methods are distinguished into two main categories. Moment matching (MM) techniques~\cite{Odabasioglu1998} are preferred due to their computational efficiency. However, they rely on an ad hoc selection of the number of moments, which correlates the final ROM size with the number of ports. On the other hand, techniques based on balanced truncation (BT)~\cite{Antoulas2004} offer reliable bounds for the approximation error and have no fundamental limitation to the number of ports they can handle, resulting in more compact ROMs. Nevertheless, BT applies only to small-scale models since it involves the computationally expensive solution of Lyapunov equations~\cite{Antoulas2004}.

In this work, appropriate performance improvements are explored to overcome the main drawback of the conventional BT method. To this end, we adopt an efficient low-rank technique based on the extended Krylov subspace (EKS) for solving the Lyapunov equations. The proposed approach can be integrated into industrial extraction tools, such as the ANSYS RaptorX\texttrademark ~\cite{raptX}, to obtain more compact ROMs of large-scale multi-port RLCk models.
\section{Background}
\label{sec:background}

Consider the modified nodal analysis (MNA) description~\cite{MNA} of an $n$-node, $m$-branch (inductive), $p$-input, and $q$-output RLCk circuit in the time domain: 
\begin{equation} \label{mna}
    \begin{aligned}
        \begin{pmatrix}
            \mathbf{G_{n}} & \mathbf{E} \\
            \mathbf{-E}^T & \mathbf{0} 
        \end{pmatrix}
        \begin{pmatrix}
            \mathbf{v}(t) \\
            \mathbf{i}(t) 
        \end{pmatrix} +
        \begin{pmatrix}
            \mathbf{C_{n}} & \mathbf{0} \\
            \mathbf{0} & \mathbf{M} 
        \end{pmatrix}
        \begin{pmatrix}
            \dot{\mathbf{v}}(t) \\
            \dot{\mathbf{i}}(t) 
        \end{pmatrix} = 
        \begin{pmatrix}
            \mathbf{B}_1 \\
            \mathbf{0} 
        \end{pmatrix}
        \mathbf{u}(t), \quad 
        \mathbf{y}(t) = 
        \begin{pmatrix}
            \mathbf{L}_1 \quad
            \mathbf{0} 
        \end{pmatrix}
        \begin{pmatrix}
            \mathbf{v}(t) \\
            \mathbf{i}(t) 
        \end{pmatrix}
    \end{aligned}
\end{equation}
where $\mathbf{G_{n}} \in\R^{n\times n}$ (node conductance matrix),  $\mathbf{C_{n}}\in\R^{n\times n}$ (node capacitance matrix), $\mathbf{M}\in\R^{m\times m}$ (branch inductance matrix), $\mathbf{E}\in\R^{n \times m}$ (node-to-branch incidence matrix), $\mathbf{v}\in\R^{n}$ (vector of node voltages), $\mathbf{i}\in\R^{m}$ (vector of inductive branch currents), $\mathbf{u} \in \R^{p} $ (vector of input excitations), $\mathbf{B}_1\in\R^{n\times p}$ (input-to-node connectivity matrix), $\mathbf{y}\in\R^{q}$ (vector of output measurements), and $\mathbf{L}_1\in\R^{q\times n}$ (node-to-output connectivity matrix). 
Moreover, we denote
$\dot{\mathbf{v}}(t) \equiv \frac{d\mathbf{v}(t)}{dt}$ and $\dot{\mathbf{i}}(t) \equiv \frac{d\mathbf{i}(t)}{dt}$.
If we now define the model order as $N \equiv n + m$, the state vector as $\mathbf{x}(t) \equiv \begin{pmatrix}
\mathbf{v}(t) \\
\mathbf{i}(t) 
\end{pmatrix}$, and~also:
\begin{equation*}
\begin{aligned}
\mathbf{G}\equiv -\begin{pmatrix}
\mathbf{G_{n}} & \mathbf{E} \\
\mathbf{-E}^T & \mathbf{0} 
\end{pmatrix},\quad  \mathbf{C} \equiv \begin{pmatrix}
\mathbf{C_{n}} & \mathbf{0} \\
\mathbf{0} & \mathbf{M} 
\end{pmatrix},\quad \mathbf{B}\equiv \begin{pmatrix}
\mathbf{B}_1 \\
\mathbf{0} 
\end{pmatrix}, \quad \mathbf{L}\equiv \begin{pmatrix}
\mathbf{L}_1 \quad
\mathbf{0} 
\end{pmatrix}
\end{aligned}, 
\end{equation*}
then Eq.~(\ref{mna}) can be written in the generalized state-space form, or so-called descriptor form:
\begin{equation}
\begin{aligned} \label{state}
\mathbf{C}\frac{d \mathbf{x}(t)}{d t} = \mathbf{G x}(t) + \mathbf{Bu}(t), \quad 
\mathbf{y}(t) = \mathbf{L x}(t).
\end{aligned}
\end{equation}

\noindent The objective of MOR is to produce an equivalent ROM: 
\begin{equation}
\begin{aligned}
\mathbf{\tilde C} \frac{d \mathbf{\tilde x}(t)}{d t} =\mathbf{\tilde G} \mathbf{\tilde x}(t) + \mathbf{\tilde B} \mathbf{u(t)}, \quad
\mathbf{\tilde y}(t) = \mathbf{\tilde L \tilde x}(t)
\end{aligned}
\end{equation}				
where  $\mathbf{\tilde G}, \mathbf{\tilde C} \in \R^{r\times r} $, $\mathbf{\tilde B} \in \R^{r\times p} $, $\mathbf{\tilde L} \in \R^{q\times r}$, the reduced order $r<<N$, and the output error is bounded as $||\mathbf{\tilde{y} }(t) -\mathbf{y}(t)||_2 < \varepsilon||\mathbf{u}(t)||_2$ for given  $\mathbf{u}(t)$ and small $\varepsilon$. The output error bound can be expressed in the frequency domain as $||\mathbf{\tilde{y} }(s) -\mathbf{y}(s)||_2 < \varepsilon||\mathbf{u}(s)||_2$ via Plancherel's theorem~\cite{Plancherel}.~If
\begin{equation*}
\begin{aligned}
\mathbf{H}(s) = \mathbf{L}(s\mathbf{C} - \mathbf{G})^{-1} \mathbf{B}, \quad
\mathbf{\tilde H}(s) =  \mathbf{\tilde L}(s\mathbf{\tilde C} - \mathbf{\tilde G})^{-1} \mathbf{\tilde B}
\end{aligned}
\end{equation*}
are the transfer functions of the original model and the ROM, the corresponding
output
error~is: 
\begin{equation}
\begin{aligned}
||\mathbf{\tilde{y} }(s) -\mathbf{y}(s)||_2 = ||\mathbf{\tilde{H}}(s) \mathbf{u}(s) - \mathbf{H}(s)\mathbf{u}(s)||_2 \quad \leq \quad ||\mathbf{\tilde{H}}(s) - \mathbf{H}(s)||_\infty||\mathbf{u}(s)||_2	
\end{aligned}
\end{equation}
where $||.||_\infty$ is the 
$\mathcal{L}_2$ matrix norm or $\mathcal{H}_\infty$ norm of a rational transfer function. 
Thus, to~bound this error, we need to bound the distance between the transfer functions:~$||\mathbf{\tilde{H}}(s) - \mathbf{H}(s)||_\infty~<~\varepsilon$.
\section{MOR by Balanced Truncation}
\label{sec:bt}

BT relies on the computation of the controllability Gramian $\mathbf{P}$ and observability Gramian~$\mathbf{Q}$, which are calculated as the solutions of the following Lyapunov matrix equations \cite{Antoulas2004}:
\begin{equation}
\begin{aligned}\label{Eq:lyap_sta}
(\mathbf{C}^{-1}\mathbf{G}) \mathbf{P} +  \mathbf{P} (\mathbf{C}^{-1}\mathbf{G})^T  = - (\mathbf{C}^{-1}\mathbf{B}) (\mathbf{C^{-1}}\mathbf{B})^T, \quad
(\mathbf{C}^{-1}\mathbf{G})^T \mathbf{Q} +  \mathbf{Q}(\mathbf{C}^{-1}\mathbf{G}) = - \mathbf{L}^T \mathbf{L}.
\end{aligned}
\end{equation}

The controllability Gramian $\mathbf{P}$ characterizes the input-to-state behavior, i.e., the degree to which the states are controllable by the inputs, while the observability Gramian $\mathbf{Q}$ characterizes the state-to-output behavior, i.e., the degree to which the states are observable at the outputs.~In principle, a ROM can be obtained by eliminating the states that are difficult to reach or observe. However, in the original state-space coordinates, there are states that are difficult to reach but easy to observe, and vice versa. The process of ``balancing'' transforms the state~vector to a new coordinate system, where for each state, the degree of difficulty is the same for both reaching and observing it. An appropriate transformation $\mathbf{Tx}(t)$ exists, leading to the following model:
\begin{equation}\label{Eq:Trunc}
\begin{aligned}
\mathbf{TCT}^{-1}\frac{d (\mathbf{Tx}(t))}{d t} = \mathbf{TGT^{-1}} (\mathbf{Tx}(t)) + \mathbf{TBu}(t), \quad
\mathbf{y}(t) = \mathbf{LT^{-1}} (\mathbf{Tx}(t))
\end{aligned}
\end{equation}
that preserves the transfer function H(s). 
This  renders
$\mathbf{P}$\ =\ $\mathbf{Q}$\ =\ $diag(\sigma_1, \sigma_2, \dots, \sigma_N)$~\cite{Antoulas2004},
~where $\sigma_i$ are known as the Hankel singular values (HSVs) of the model and are equal to the square roots of the eigenvalues of product $\mathbf{P Q}$,  
i.e., $\sigma_i = \sqrt{\lambda_i(\mathbf{P Q})}$.
In the above balanced model,
the states that are easier to reach and observe correspond to the largest HSVs. If $r$ of them are preserved (truncating the $N-r$ states corresponding to the smallest HSVs), it can be shown that the distance between the original and the reduced-order transfer functions is bounded as:
\begin{equation}\label{Eq:errr}
||\mathbf{H}(s) - \mathbf{\tilde H}(s)||_\infty \leq 2( \sigma_{r+1} +\sigma_{r+2}+...+\sigma_{N})
\end{equation}
The latter is an ``a-priori'' criterion for selecting the ROM order for a desired output error tolerance $\varepsilon$, which constitutes a significant advantage of BT over MM techniques. 
The main steps of the BT procedure are summarized in Algorithm \ref{vanilla_BT}.
\begin{figure}[!hbt]
\vspace{-1.5em}
\begin{algorithm}[H]
\footnotesize
	\caption{MOR by balanced truncation}\label{vanilla_BT}
	\begin{algorithmic}[1]
		\State Solve the Lyapunov equations to obtain the Gramian matrices $\mathbf{P}$ and $\mathbf{Q}$~\cite{Lathauwer2004}
		\State Compute the
            SVD of the Gramian matrices: $\mathbf{P} = \mathbf{U}_P \mathbf{\Sigma}_P \mathbf{V}_P^T $ and $\mathbf{Q} = \mathbf{U}_Q \mathbf{\Sigma}_Q \mathbf{V}_Q^T $
            \State Find the square root of the Gramian matrices: $\mathbf{Z}_P = \mathbf{U}_P \mathbf{\Sigma}_P^{1/2}$ and $\mathbf{Z}_Q = \mathbf{U}_Q \mathbf{\Sigma}_Q^{1/2}$
            \State Compute the SVD of the product of the roots: $\mathbf{Z}_Q^T\mathbf{Z}_P = \mathbf{U}\mathbf{\Sigma}\mathbf{V}^T$	
            \State Compute transformation matrices: $\mathbf{T}_{(r\times N)}$\ =\ $\mathbf{\Sigma}_{(r\times r)}^{-1/2} \mathbf{U}_{(r\times N)} \mathbf{Z}_Q^T $,\ $\mathbf{T}_{(N\times r)}^{-1}$\ =\ $ \mathbf{Z}_P\mathbf{V}_{(N\times r)}\mathbf{\Sigma}^{-1/2}_{(r\times r)}$
            \State Compute ROM: $\mathbf{\tilde G}$\ =\ $\mathbf{T}_{(r\times N)}\mathbf{G}\mathbf{T}_{(N\times r)}^{-1}$, \ $\mathbf{\tilde C}$\ =\ $\mathbf{T}_{(r\times N)}\mathbf{C}\mathbf{T}_{(N\times r)}^{-1}$,\  $\mathbf{\tilde B}$\ =\ $\mathbf{T}_{(r\times N)}\mathbf{B},\ \mathbf{\tilde L}$\ =\ $\mathbf{L}\mathbf{T}_{(N\times r)}^{-1}$
	\end{algorithmic}
\end{algorithm}
\vspace{-1.5em}
\end{figure}
The major drawback of BT is the significant computational and memory cost for deriving the ROM, which hinders the applicability to large-scale models (with $N$ over a few thousand states). 
This is because the operations involved (e.g., the solution of Lyapunov equations and the singular value decomposition [SVD]) are computationally expensive with a complexity of $O(N^3)$. Moreover, they are applied on dense matrices, since the Gramians $\mathbf{P}, \mathbf{Q}$ are dense even if the system matrices $\mathbf{C}, \mathbf{G}, \mathbf{B}, \mathbf{L}$~are~sparse. 

However, the products $(\mathbf{C}^{-1}\mathbf{B}) (\mathbf{C^{-1}}\mathbf{B})^T$ and $\mathbf{L}^T\mathbf{L}$ have low numerical order compared~to~$N$, as $p,q<<N$, resulting in low-rank Gramian matrices that can be approximated using low-rank techniques. This greatly reduces the complexity and memory requirements of the solution of~the Lyapunov equations and the SVD analysis, which are now of order $k$ instead of full order~$N$.

\subsection{Low-rank BT MOR}
\label{sec:bt_eks}
The essence of low-rank BT MOR is to iteratively project the Lyapunov equations of Eq.~(\ref{Eq:lyap_sta}) onto a lower-dimensional Krylov subspace and then solve the resulting small-scale equations to obtain low-rank approximate solutions of Eq.~(\ref{Eq:lyap_sta}). More specifically, if $\mathbf{K} \in \R^{N \times k}$ ($k<<N$) is a projection matrix whose columns span the $k$-dimensional Krylov subspace:
 \begin{equation*} 
\mathcal{K}_k(\mathbf{G}_{C},\mathbf{B}_{C})=  
  span  \{\mathbf{B}_{C},\mathbf{G}_{C}\mathbf{B}_{C},  \mathbf{G}_{C}^{2}\mathbf{B}_{C},\dots,\mathbf{G}_{C}^{k-1}\mathbf{B}_{C}\}
 \end{equation*}
where $\mathbf{G}_{C} \equiv \mathbf{C}^{-1}\mathbf{G},\hspace{0.25em} \mathbf{B}_{C} \equiv \mathbf{C}^{-1}\mathbf{B}$,
then the projected Lyapunov equation (for the controllability Gramian $\mathbf{P}$) onto $\mathcal{K}_k(\mathbf{G}_{C},\mathbf{B}_{C})$ is:
\begin{equation} \label{Eq:projected_lyap}
(\mathbf{K}^T\mathbf{G}_{C}\mathbf{K})\mathbf{X} +\mathbf{X} (\mathbf{K}^T\mathbf{G}_{C}\mathbf{K})^T   =-\mathbf{K}^T\mathbf{B}_{C}\mathbf{B}_{C}^T\mathbf{K}
\end{equation}
(the same holds true for the observability Gramian $\mathbf{Q}$ with $\mathbf{G}_{C}^T$, $\mathbf{L}^{T}$ in place of $\mathbf{G}_{C}$, $\mathbf{B}_{C}$). The solution $\mathbf{X} \in \R^{k \times k}$ of  Eq.~(\ref{Eq:projected_lyap}) can be back-projected to the $N$-dimensional space to give an approximate solution $\mathbf{P} = \mathbf{K}\mathbf{XK}^T$ for the original large-scale Eq.~(\ref{Eq:lyap_sta}), and a  low-rank factor $\mathbf{Z} \in \R^{N \times k}$ of $\mathbf{P}$ can be obtained as $\mathbf{Z} =\mathbf{K}\mathbf{U}\mathbf{\Sigma}^{1/2}$, where $[\mathbf{U},\mathbf{\Sigma},\mathbf{V}] = SVD(\mathbf{X})$. 

Although the projection process is independent of the subspace selection, its effectiveness is critically dependent on the chosen subspace.
The convergence to the final solution
can be accelerated by enriching the standard Krylov subspace $\mathcal{K}_k(\mathbf{G}_{C},\mathbf{B}_{C}) $ with information from the subspace $\mathcal{K}_k(\mathbf{G}_{C}^{-1},\mathbf{B}_{C})$, which corresponds to the inverse matrix $\mathbf{G}_{C}^{-1}$, leading to the EKS~\cite{Arxiv_22, ASPDAC21}:
\begin{equation} \label{Eq:eksm}
\mathcal{K}_k^C(\mathbf{G}_{C},\mathbf{B}_{C}) = 
span \{\mathbf{B}_{C},  \mathbf{G}_{C}^{-1}\mathbf{B}_{C}, \mathbf{G}_{C}\mathbf{B}_{C},\mathbf{G}_{C}^{-2}\mathbf{B}_{C}, \mathbf{G}_{C}^{2}\mathbf{B}_{C},\dots, \\
\mathbf{G}_{C}^{-(k-1)}\mathbf{B}_{C}, \mathbf{G}_{C}^{k-1}\mathbf{B}_{C}\}
\end{equation}
The EKS method (EKSM) starts with the pair $\{ \mathbf{B}_{C}, \mathbf{G}_{C}^{-1}\mathbf{B}_{C} \}$ and generates an extended subspace $\mathcal{K}_k^C(\mathbf{G}_{C},\mathbf{B}_{C})$ of increasing dimension, solving the projected Lyapunov Eq.~(\ref{Eq:projected_lyap}) in each iteration, until a sufficiently accurate approximation of the solution of Eq.~(\ref{Eq:lyap_sta}) is obtained.~The complete EKSM is presented in Algorithm~\ref{eksm-alg}. 
Below are some efficient implementation details:
\begin{itemize}
\itemsep0em
    \item \textbf{Matrix inversion by linear solves}: 
	 The inputs to Algorithm \ref{eksm-alg} are not actually $\mathbf{G}_{C} \equiv \mathbf{C}^{-1}\mathbf{G}$ or $\mathbf{G}^T_{C} \equiv (\mathbf{C}^{-1}\mathbf{G})^T$ but the system matrices $\mathbf{G}$, $\mathbf{C}$ or $\mathbf{G}^T$, $\mathbf{C}^T$, since the (generally dense) inverse matrices are only needed in products with $p$ vectors (in step ~\hyperref[eksm-alg:step2]{2}) and $2pj$ vectors (in steps ~\hyperref[eksm-alg:step4]{4} and ~\hyperref[eksm-alg:step11]{11} of each iteration).
  These can be implemented as  linear solves $\mathbf{C}\mathbf{Y} = \mathbf{R}$ and $\mathbf{G}\mathbf{Y} = \mathbf{R}$ (or $\mathbf{C}^T\mathbf{Y} = \mathbf{R}$, $\mathbf{G}^T\mathbf{Y} = \mathbf{R}$) by any direct or iterative algorithm like~\cite{Bavier2012}. 
  \item \textbf{Handling of sparse/dense matrices}: Note that matrix $\mathbf{M}$ of Eq.~(\ref{mna}) is highly dense,~as it generally includes a huge number of mutual inductances. To effectively handle the sparse ($\mathbf{C}_n$) and dense ($\mathbf{M}$) blocks of matrix $\mathbf{C}$, we use 
  efficient data structures and numerical~techniques. 
  For example,  for linear solves and matrix-vector products, we employ parallel CPU-optimized methods for sparse matrices and leverage GPU-accelerated techniques~\cite{CMG_GPU} for dense matrices.
  \item \textbf{Solution of the small-scale Lyapunov equations}: 
  To solve the small-scale~($2pj \times 2pj$) Lyapunov~equations  in step~\hyperref[eksm-alg:step5]{5} of each iteration,
  we employ the Bartels-Stewart~algorithm~\cite{Lathauwer2004}. 
  \item  \textbf{Convergence criterion:} An appropriate stopping criterion is the residual of Eq.~(\ref{Eq:lyap_sta}) with the approximate solution $\mathbf{P} = \mathbf{K}\mathbf{XK}^T$ to reach a certain threshold in magnitude, i.e.,
	\begin{equation}\label{Eq:crit}
	\frac{||\mathbf{G}_{C}\mathbf{K}^{(j)}\mathbf{X}\mathbf{K}^{(j)T} + \mathbf{K}^{(j)}\mathbf{X}\mathbf{K}^{(j)T}\mathbf{G}_{C} + \mathbf{B}_{C}\mathbf{B}_{C}^{T}||}{||\mathbf{B}_{C}\mathbf{B}_{C}^{T}||} \le tol 
	\end{equation}
However, this criterion equals to $||\mathbf{R}^T\mathbf{M}\mathbf{X}|| \le tol$~\cite{Simoncini2007}, which can be computed more efficiently. A tolerance of $tol$\ =\ $10^{-10}$ is typically adequate to obtain an accurate model.
\end{itemize}
\begin{figure}[!hbt]
\vspace{-1.5em}
\begin{algorithm}[H]
\footnotesize
\caption{Extended Krylov subspace method for low-rank solution of Lyapunov equations}\label{eksm-alg}	
	\textbf{Input:}  $ \mathbf{G}_{C} \equiv \mathbf{C}^{-1}\mathbf{G},  \mathbf{B}_{C} \equiv \mathbf{C}^{-1}\mathbf{B}$  (or  $\mathbf{G}_{C}^T$, $\mathbf{L}^T $)\\
	\textbf{Output:} $\mathbf{Z}$ such that $\mathbf{P} \approx \mathbf{Z} \mathbf{Z} ^T $
	\begin{algorithmic} [1]
		\State $j=1$; $p=size\_col(\mathbf{B}_{C})$
        \label{eksm-alg:step2}
		\State $\mathbf{K}^{(j)} = Orth([\mathbf{B}_{C},\mathbf{G}_{C}^{-1}\mathbf{B}_{C}])$
		\While {$j < maxiter$}
        \label{eksm-alg:step4}
    		\State $\mathbf{A} = \mathbf{K}^{(j)T}\mathbf{G}_{C}\mathbf{K}^{(j)} $;\quad $\mathbf{R} = \mathbf{K}^{(j)T}\mathbf{B}_{C} $  
        \label{eksm-alg:step5}
    		\State Solve $\mathbf{A}\mathbf{X} +\mathbf{X} \mathbf{A}^{T}  = -\mathbf{R}\mathbf{R}^{T}$ for  $\mathbf{X} \in \R^{2pj \times 2pj}$
    		\If {converged}
        		\State $[\mathbf{U}, \mathbf{\Sigma}, \mathbf{V}] = \mathbf{SVD}(\mathbf{X})$; \quad
        		$\mathbf{Z} = \mathbf{K}^{(j)}\mathbf{U}\mathbf{\Sigma}^{1/2}$
        	    \State \textbf{break}
        	\EndIf
    		\State $k_1=2p(j-1)$; $k_2 = k_1+p$; $k_3 = 2pj$ 
        \label{eksm-alg:step11}
    		\State $\mathbf{K}_1  = [\mathbf{G}_{C}\mathbf{K}^{(j)}(:,k_1+1:k_2),\mathbf{G}_{C}^{-1}\mathbf{K}^{(j)}(:,k_2+1:k_3)]$
    		\State $\mathbf{K}_2 = Orth(\mathbf{K}_1) $ \quad w.r.t. \quad $\mathbf{K}^{(j)}$
    		\State $\mathbf{K}_3 = Orth(\mathbf{K}_2) $ 
    		\State $\mathbf{K}^{(j+1)} = [\mathbf{K}^{(j)},\mathbf{K}_3]$\quad 
		    \State $j=j+1$
		\EndWhile
	\end{algorithmic}
\end{algorithm}
\vspace{-1.5em}
\end{figure}

\section{Experimental Evaluation}
\label{sec:experimental}
\subsection{Experimental setup}
\label{sec:exp_setup}

To evaluate EKSM, we used large-scale RLCk models extracted from different circuits~using ANSYS RaptorX\texttrademark~\cite{raptX}. These circuits consist of many passive elements, including mutual~inductances. The EKSM ROMs are compared against golden ROMs produced by~RaptorX\texttrademark, through S-parameter plotting. The characteristics of the RLCk models are listed in~Table~\ref{table:AUTH_models}. All experiments were executed on a Linux server with a 2.80 GHz 16-thread CPU and 64~GB~of~memory.

\begin{table}[!hbt]
\centering
        \vspace{-0.4em}
	\caption{ Detailed characteristics of RLCk models }
        \vspace{-7pt}	
 \label{table:AUTH_models}
\small
\setlength{\tabcolsep}{3pt}
\begin{tabular}{|c|c|c|c|c|c|c|c|}
\hline
 Model & Initial order & \#nodes & \#ports & \#resistors & \#capacitors & \#inductors & \#mutual ind. \\ \hline
 VGA\_28 & 95189 & 57675 & 13 & 155879 & 169600 & 37514 & 126766838 \\ \hline
 Hybrid\_56 & 98024 & 59210 & 5 & 112338 & 290572 & 38814 & 165802476 \\ \hline
 Wilkinson\_56 & 100888 & 60703 & 4 & 115117 & 271293 & 40185 & 193641938 \\ \hline
 VCO\_13 & 104367 & 61264 & 4 & 604072 & 596846 & 43103 & 188436057 \\ \hline
 CSLNA\_56 & 128574 & 78046 & 9 & 188842 & 472573 & 50528 & 169339965 \\ \hline
 Wilkinson\_28 & 129087 & 78263 & 4 & 123254 & 266710 & 50824 & 259462454 \\ \hline
 Hybrid\_28 & 134710 & 75766 & 5 & 128935 & 283905 & 53169 & 264162513 \\ \hline
 LNACASC\_28 & 162881 & 96876 & 11 & 774427 & 684662 & 66005 & 323090671 \\ \hline
\end{tabular}
\vspace{-1.25em}
\end{table}
\subsection{Experimental results}
\label{sec:results_bt_low-rank}

\begin{figure}[!tbh]
  \centering
    \vspace{-1em}
    \includegraphics[width=0.87\textwidth]{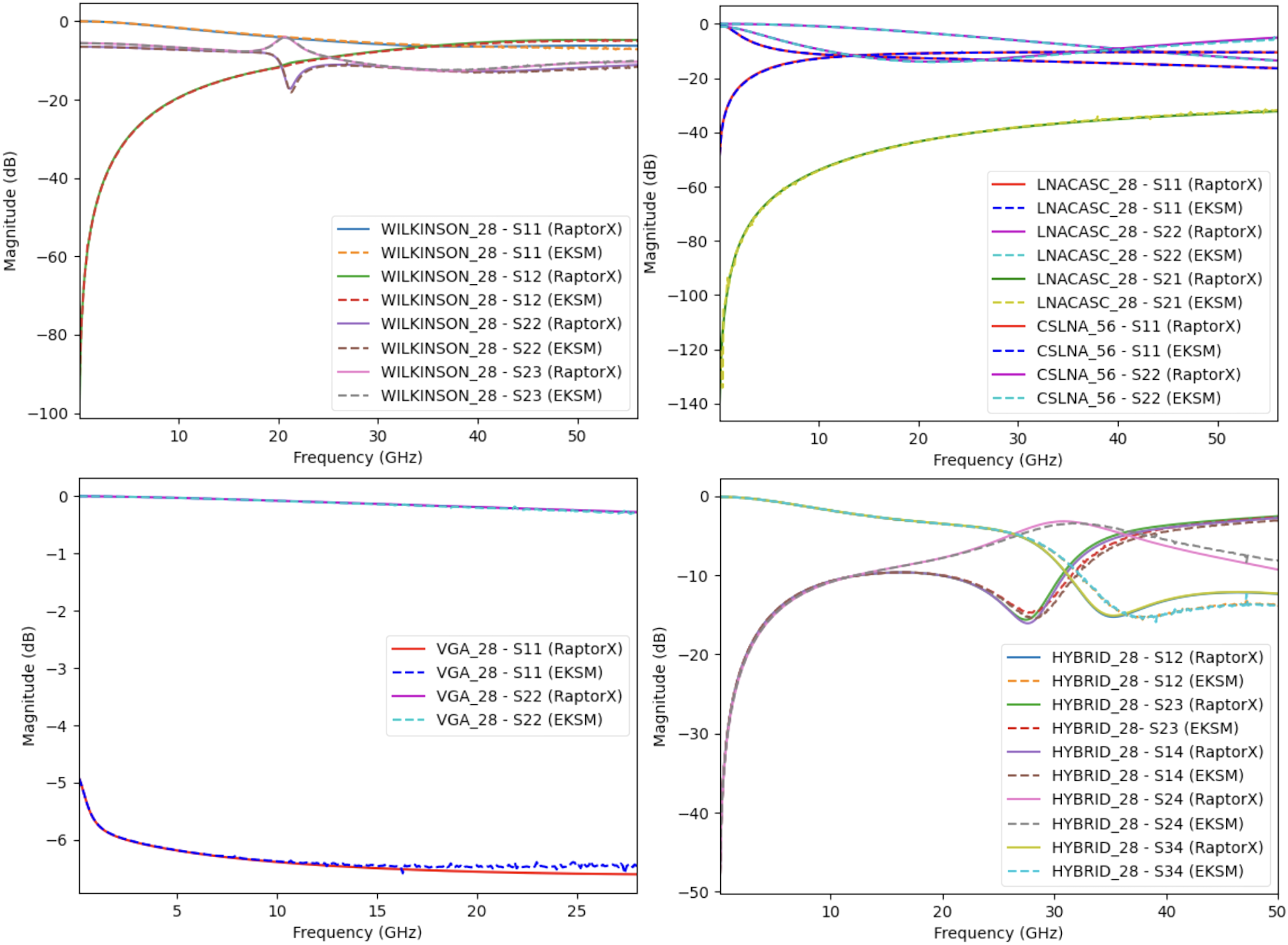}
    \vspace{-8pt}
    \caption{Comparison of accuracy between EKSM and RaptorX\texttrademark~ROMs.}
    \label{fig:accuracy_results}
    \vspace{-1.8em}
\end{figure}

The efficiency of the EKSM against RaptorX\texttrademark\ is demonstrated in Table~\ref{table:low_rank_results}. The S-parameters plots of Figure~\ref{fig:accuracy_results} indicate that EKSM achieves accuracy close to that of RaptorX\texttrademark\ while producing roughly $\times$3.1 more compact ROMs. Although EKSM has higher reduction~time and memory requirements, they are still reasonable and can be significantly improved in future~work.
\begin{table}[!hbt]
\centering
    \vspace{-0.5em}
    \caption{ ROM order and MOR performance of EKSM vs RaptorX\texttrademark}
      \vspace{-6pt}
    \label{table:low_rank_results}
\small
\begin{tabular}{|c||c||c|c||c|c||c|c|}
\hline
\multirow{2}{*}{Model}  &  \multirow{2}{*}{\makecell{Initial \\ order}} & \multicolumn{2}{c||}{ROM order} & \multicolumn{2}{c||}{Reduction time (s)} & \multicolumn{2}{c|}{Memory (GB)} \\ 
\cline{3-8}& & RaptorX\texttrademark& EKSM & RaptorX\texttrademark& EKSM & RaptorX\texttrademark& EKSM\\ \hline
VGA\_28 & 95189 & 4744 & 1040 & 67 & 1037 & 32.63 & 19.14 \\ \hline
Hybrid\_56 & 98024 & 1267 & 397 & 104 & 613 & 24.05 & 29.11 \\ \hline
Wilkinson\_56 & 100888 & 765 & 320 & 154 & 570 & 24.79 & 29.76 \\ \hline
VCO\_13 & 104367 & 407 & 311 & 119 & 673 & 26.48 & 29.18 \\ \hline
LNACS\_56 & 128574 & 2172 & 716 & 74 & 1237 & 25.82 & 26.74 \\ \hline
Wilkinson\_28 & 129087 & 885 & 302 & 205 & 801 & 25.35 & 36.21 \\ \hline
Hybrid\_28 & 134710 & 787 & 399 & 217 & 1032 & 24.31 & 35.52 \\ \hline
LNACasc\_28 & 162881 & 4768 & 879 & 373 & 2866 & 78.52 & 48.67 \\ \hline
\end{tabular}
\vspace{-0.6em}
\end{table}

\section{Conclusions}
\label{sec:conclusions}
\vspace{-0.4em}
Alternative MOR techniques to reduce large-scale RLCk models with accuracy comparable~to commercial tools are presented. The proposed low-rank BT method is evaluated across diverse large-scale benchmark circuits by comparing their S-parameters. Experimental results indicate that our approach achieves sufficient accuracy~while providing ROMs that are up to ×5.5 smaller than the ROMs obtained by~ANSYS~RaptorX\texttrademark.
\vspace{-0.6em}
\section{Acknowledgments}
\label{sec:acknowledgments}
\vspace{-0.4em}
This research has been co-financed by the European Regional Development Fund and Greek national funds via the Operational Program "Competitiveness, Entrepreneurship and Innovation,"~under the call "RESEARCH-CREATE-INNOVATE" (project code: T2EDK-00609).
\vspace{-1.6em}

\bibliographystyle{IEEEtran}
\bibliography{easychair}

\end{document}